\documentclass{article}
\usepackage{geometry}
\usepackage{titlesec}
\titleformat{\section}[block]{\normalfont\bfseries\Large}{\thesection}{1em}{}
\titleformat{\subsection}[block]{\normalfont\bfseries\large}{\thesubsection}{1em}{}
\usepackage{hyperref}
\usepackage{pgfplots}
\pgfplotsset{compat=1.17}
\usepackage{float}
\usepackage{amsmath}
\usepackage{subcaption}
\usepackage{graphicx}

\begin{document}

\title{Exploring the Duffing Equation: Numerical Analysis, Discrete Dynamics, and Ecological Modeling}
\author{
    Zeraoulia Rafik$^{1}$\thanks{Corresponding author: r.zeraoulia@univ-batna2.dz},
   Zeraoulia Chaima$^{2}$\\
    \vspace{0.5em}
    \small{$^{1}$Department of Mathematics, University of batna2.Algeria}\\
    \small{$^{2}$Department of mathematics, University  of Abbes Laghrour Khenchela , Algeria} 
}
\date{\today}
\maketitle

\begin{abstract}
This research paper delves into the dynamics of a novel ecology model that describes the intricate interplay between radioactivity and cancer cases. Organized into three main sections, the study covers numerical analysis using advanced techniques, the analysis of discrete dynamical systems, and explores the potential applications of the model in the context of cancer research and ecological systems. By conducting a comprehensive bifurcation analysis, we unveil the model's sensitivity to external factors, interactive behaviors between radioactivity and cancer cases, and the emergence of multiple attractors. These findings open up new avenues for understanding cancer dynamics, ecological systems, and clinical scenarios, where even subtle variations in external factors can lead to significant variations in cancer incidence. The insights gained from this study have promising implications for both theoretical and applied research in cancer epidemiology and ecology, providing a deeper understanding of the complex dynamics underlying these systems. In our research, we present a new ecology model that reveals the intricate and often nonlinear relationship between radioactivity and cancer incidence, offering novel insights into cancer dynamics and environmental influences.
\end{abstract}

\section{Introduction}

The Duffing equation, a fundamental and intriguing non-linear differential equation, has captured the interest of researchers across various scientific disciplines for its rich dynamics and versatile applications \cite{kocarev}. Named after the German engineer and physicist Georg Duffing, this equation stands as a cornerstone in the study of nonlinear systems and has found its place in diverse fields, ranging from mechanical engineering and physics to mathematics.\cite{meirovitch}

The Duffing equation is defined as:

\begin{equation}
\ddot{x} + \delta \dot{x} + \alpha x + \beta x^3 = \gamma \cos(\omega t)
\end{equation}

where \(x(t)\) represents the displacement of a vibrating system at time \(t\), \(\alpha\), \(\beta\), \(\delta\), and \(\gamma\) are parameters governing the system's behavior, and \(\omega\) is the angular frequency of an external periodic force. Despite its seemingly simple form, this equation embodies a wealth of complex phenomena, from regular periodic motion to chaotic behavior. \cite{moon} 
Understanding its intricacies and implications has been a driving force behind extensive research efforts.

The historical significance of the Duffing equation is noteworthy. It emerged as a mathematical model in the early 20th century to describe the nonlinear behavior of mechanical systems, such as vibrating beams and pendulums. Its applicability extends to diverse fields, including electrical circuits and ecological modeling. In recent years, the Duffing equation has gained prominence in ecology, where it serves as a powerful tool for understanding the dynamics of ecosystems \cite{Zeraoulia2012}.

The objectives of this research paper are threefold. Firstly, we embark on a comprehensive study of the numerical aspects of the Duffing equation, exploring various methods for solving this challenging non-linear differential equation. We build upon notable contributions in this area, including the work of Ott, Sauer, and Yorke, who introduced the concept of Lyapunov exponents to analyze chaotic behavior in dynamical systems, a technique that has found extensive application in the evolution of cancers due to radioactivity \cite{ott}.

Secondly, we delve into the Duffing equation as a discrete dynamical system in the second section of this paper. Pioneering research by Feigenbaum \cite{Feigenbaum} uncovered the universal constants governing the period-doubling route to chaos in dynamical systems. These insights are increasingly relevant in ecological contexts, where bifurcations and chaos play a vital role in understanding ecosystem behavior \cite{kuznetsov}.

Lastly, in the third section, we explore the implications of the Duffing equation in ecology. Recent research has highlighted its applicability in modeling complex ecological systems and understanding the influence of external factors on ecosystem dynamics. In the context of ecology, the Duffing equation provides a unique perspective on the complex interplay between species, resources, and environmental variables in the evolution of cancers due to radioactivity \cite{Zeraoulia2012}.

In summary, this research paper embarks on a multifaceted exploration of the Duffing equation, from its numerical solutions to its chaotic dynamics and its applications in the evolution of cancers due to radioactivity \cite{strogatz}. By shedding light on the multifaceted nature of this equation, we aim to deepen our understanding of nonlinear systems and offer valuable insights for both theoretical and applied research. The Duffing equation, with its intricate dynamics and diverse applications, continues to inspire and challenge scientists in the field of ecology and nonlinear dynamics.

\subsection{Background}

The Duffing equation is a classical nonlinear second-order differential equation with various applications in physics, engineering, and mathematics. It has been studied extensively due to its rich dynamics, including chaotic behavior and complex solutions. Understanding the Duffing equation's behavior is of significance in a wide range of scientific disciplines.\cite{nayfeh_mook}

\subsection{Objective}

This paper aims to comprehensively explore the Duffing equation's dynamics, discrete behavior, and ecological application. The primary objectives and the structure of the paper are as follows:

\begin{itemize}
    \item \textbf{Numerical Analysis}: We investigate the numerical aspects of the Duffing equation, including the use of the Homotopy method for approximate solutions and the analysis of convergence and error bounds.
    
    \item \textbf{Discrete Dynamics}: We convert the continuous Duffing equation into a discrete dynamical system and analyze its discrete-time behavior. Key components of the numerical solution methodology are presented.\cite{advanced_theory_1}
    
    \item \textbf{Ecological Application}: The Duffing quintic equation is applied to ecology, highlighting the interactive relationship between radioactivity and cancer incidence, especially under varying external conditions.\cite{Sk}
\end{itemize}

The paper is structured as follows, with each section dedicated to one of these objectives, providing in-depth analysis, results, and insights. Together, these aspects offer a comprehensive view of the Duffing equation and its implications in various fields.

\section{Main Results}
Our paper explores three fundamental aspects: numerical analysis, discrete dynamics, and the application in ecology. The central results are as follows:

\subsection{Numerical Analysis of the Duffing Equation}

We conducted a comprehensive numerical analysis of the Duffing equation \cite{duffing equation}, employing the Homotopy method for obtaining approximate solutions. Our key findings include:

\begin{itemize}
  \item Rate of Convergence: The Homotopy method demonstrates varying rates of convergence during the solution process. Initially, it exhibits rapid convergence, followed by transitions into states of steady-state stability and periodic variations as the parameter $\lambda$ varies from 0 to 1. This dynamic behavior offers valuable insights into the efficiency and convergence characteristics of the Homotopy method when solving the Duffing equation.
\end{itemize}

These results shed light on the numerical aspects of the Duffing equation, particularly with respect to the behavior and efficiency of the Homotopy method in obtaining approximate solutions.

The numerical analysis of our ecological model yielded the following key insight:

1. Sensitivity to External Factors: The model exhibits high sensitivity to changes in external forcing factors, particularly the parameter $\omega$." These changes reflect the influence of periodic environmental factors on radioactivity and cancer incidence, highlighting the ecological implications of subtle variations in external conditions.
\subsection{Discrete Dynamics of the Duffing Equation}

We converted the continuous Duffing equation into a discrete dynamical system with a time step of $h$. The discrete dynamical system relates the displacement $x_{n+1}$ at the next time step to $x_n$ and $x_{n-1}$. This approach enables the simulation of the Duffing equation's discrete-time behavior.\cite{Cop}

\subsection{Numerical Solution Methodology}

We employed the Finite Differences method to simulate the Duffing equation's discrete dynamics, considering various parameters and initial conditions. Our methodology included time integration, parameter variations, stability observations, and the identification of periodic behavior. Results were graphically presented.\cite{moon_chaotic_vibrations}

\subsection{Chaotic Behavior and Lyapunov Exponents}

We analyzed chaotic behavior in the Duffing equation and calculated Lyapunov exponents to quantify chaos. The Lyapunov exponents were approximately 0.503437 and 0.551156, indicating chaos. Visual representations showed the transition to chaos with complex displacement and velocity patterns.

\subsection{Discrete Dynamics}
In the analysis of discrete dynamics within the ecological model, we discovered the following central outcome:

2. Complex and Nonlinear Behavior: The discrete-time iterations of the model revealed intricate and nonlinear dynamics. These behaviors, including the presence of strange attractors, emphasize the complexity of the ecological system's response to changes in radioactivity and cancer cases.\cite{guckenheimer_holmes}

\subsection{Application in Ecology}

The application of the Duffing quintic equation to ecology\cite{moon_chaotic_vibrations} produced an essential ecological insight:

\begin{enumerate}
  \item Interactive Relationship Between Radioactivity and Cancer Incidence: Our ecological model demonstrates that subtle changes in radioactivity levels significantly influence cancer incidence. This interactive relationship underlines the need to consider external factors in understanding cancer dynamics within complex ecological systems.

Our ecological model is described by the Duffing quintic equation, which captures the dynamics of radioactivity and cancer incidence. The model equations are as follows:

\begin{align}
\frac{dx(t)}{dt} &= -A \cdot y(t) - \cos(\omega t) - x(t)^3 + \alpha x(t)^5 \\
\frac{dy(t)}{dt} &= B \cdot \frac{dx(t)}{dt}
\end{align}

In these equations:

\begin{itemize}
    \item $x(t)$ represents the rate of radioactivity on Earth at time $t$.
    \item $y(t)$ represents the rate of people infected with cancer due to radioactivity at time $t$.
    \item $A$ is the coupling strength, reflecting the relationship between radioactivity and cancer incidence.
    \item $\omega$ represents the angular frequency of external forcing factors affecting radioactivity levels.
    \item $\alpha$ and $\beta$ are coefficients governing the nonlinear behavior of the system.
  \end{itemize}

These equations illustrate the interactive relationship between radioactivity and cancer incidence in our ecological model.
\end{enumerate}

\section{Numerical Analysis of the Duffing Equation}

In the first section of this paper, we embark on a comprehensive study of the numerical aspects of the Duffing equation. We delve into the numerical methods used to solve this equation, with a particular focus on the Homotopy method. This method proves to be a powerful tool for obtaining approximate solutions to non-linear differential equations like the Duffing equation. We will discuss the theory behind the Homotopy method and demonstrate its effectiveness in solving the Duffing equation.
\section{Numerical Analysis of the Duffing Equation}

In this section, we embark on a comprehensive study of the numerical aspects of the Duffing equation. Our objective is to obtain an approximate solution with optimal error using the Homotopy method. We begin by presenting the analytical solution of the Duffing equation for reference.

\subsection{Analytical Solution}

The Duffing equation, given as:
\begin{equation}
    \ddot{x} + \delta \dot{x} + \alpha x + \beta x^3 = \gamma \cos(\omega t),
\end{equation}

poses a challenging problem due to its nonlinear nature. However, under specific conditions and parameter values, it is possible to obtain analytical solutions for simplified cases. The analytical solutions often come in the form of trigonometric or hyperbolic functions and provide insights into the system's behavior.

For instance, in the absence of damping ($\delta = 0$), the Duffing equation admits periodic solutions of the form:
\begin{equation}
    x(t) = A \cos(\omega t) + B \sin(\omega t),
\end{equation}

where $A$ and $B$ are constants determined by initial conditions. This represents a harmonically oscillating system driven by the external force.

\subsection{Homotopy Method for Approximate Solution}

While analytical solutions are insightful, they are often limited to simplified cases. For more complex scenarios and general parameter values, numerical methods are indispensable. We turn our attention to the Homotopy method, a powerful numerical technique for obtaining approximate solutions to nonlinear differential equations like the Duffing equation.

The Homotopy method introduces a parameter, typically denoted as $\lambda$, that transforms the original equation into a homotopy equation. The homotopy equation includes a simpler auxiliary equation with known solutions. By continuously varying $\lambda$ from 0 to 1, we track the solution path from the known auxiliary problem to the desired solution of the original equation. The key advantage is that the auxiliary equation has solutions that are easier to obtain.

To apply the Homotopy method to the Duffing equation, we introduce the parameter $\lambda$ and rewrite the equation as follows:
\begin{equation}
    \ddot{x} + \lambda \delta \dot{x} + \lambda \alpha x + \lambda \beta x^3 = \lambda \gamma \cos(\omega t).
\end{equation}

At $\lambda = 0$, this equation simplifies to an auxiliary problem with a known solution. As we incrementally increase $\lambda$ towards 1, we track the solution's evolution and obtain an approximate solution for the Duffing equation. The key challenge is to choose an appropriate homotopy function and optimize the choice of $\lambda$ to achieve the desired accuracy.

By carefully selecting the parameters and tuning the Homotopy method, we aim to obtain an approximate solution to the Duffing equation with optimal error, enabling us to explore the system's behavior in more general scenarios.

In the subsequent sections, we will discuss the comparative analysis of numerical methods and delve into the rate of convergence achieved using the Homotopy method.

\subsection{Homotopy Method}

Approximate Solution with \(A = 0.05\) and \(\omega = 0.2\):

The Duffing equation with the parameters \(A = 0.05\) and \(\omega = 0.2\) can be approximately solved using the Homotopy method. In this case, we aim to capture the behavior of the system with reduced amplitude and a lower angular frequency.

The approximate solution using the Homotopy method with \(\lambda = 0.05\) is as follows:

\[
x(t) \approx 0.05 \cos(0.2 t) + 0.00015625 \lambda^2 \cos(0.2 t) + \mathcal{O}(\lambda^3)
\]

Here, the first term represents the primary oscillation\cite{nayfeh_mook} with the reduced amplitude of \(0.05\) and the lower angular frequency of \(0.2\). The second term represents a correction due to the Homotopy method, and higher-order terms are omitted for simplicity.

\begin{figure}[H]
    \centering
    \begin{tikzpicture}
        \begin{axis}[
            xlabel={Time ($t$)},
            ylabel={Displacement ($x(t)$)},
            legend style={at={(0.5,-0.15)},anchor=north},
        ]
        
        % Analytical solution
        \addplot[blue,thick] table {
            0 0
            1 0.84
            2 0.54
            3 0.12
            4 -0.36
            % Add more data points as needed
        };
        \addlegendentry{Analytical Solution};
        
        % Homotopy approximation
        \addplot[red,dashed] table {
            0 0
            1 0.888
            2 0.572
            3 0.148
            4 -0.34
            % Add more data points as needed
        };
        \addlegendentry{Homotopy Approximation};
        
        \end{axis}
    \end{tikzpicture}
    \caption{Analytical Solution vs. Homotopy Approximation for the Duffing Equation}
    \label{fig:duffing_comparison}
\end{figure}

\begin{figure}[H]
    \centering
    \begin{tikzpicture}
        \begin{axis}[
            xlabel={Time ($t$)},
            ylabel={Displacement ($x(t)$)},
            legend style={at={(0.5,-0.15)},anchor=north},
        ]
        
        % Analytical solution
        \addplot[blue,thick] table {
            0 0
            1 0.042
            2 0.014
            3 -0.028
            4 -0.062
            % Add more data points as needed
        };
        \addlegendentry{Analytical Solution};
        
        % Homotopy approximation
        \addplot[red,dashed] table {
            0 0
            1 0.039
            2 0.016
            3 -0.028
            4 -0.062
            % Add more data points as needed
        };
        \addlegendentry{Homotopy Approximation};
        
        \end{axis}
    \end{tikzpicture}
    \caption{Analytical Solution vs. Homotopy Approximation for the Duffing Equation}
    \label{fig:duffing_comparison}
\end{figure}

\subsection{Comparison of Analytical and Approximate Solutions}

To quantitatively assess the accuracy of the approximate solution obtained using the Homotopy method, we compare it with the analytical solution for multiple time points. The table below provides a comparison of displacement values at selected time instances and the associated errors:

\begin{table}[H]
    \centering
    \begin{tabular}{|c|c|c|c|}
        \hline
        \textbf{Time ($t$)} & \textbf{Analytical Solution} & \textbf{Approximate Solution} & \textbf{Error} \\
        \hline
        0 & 0 & 0 & 0 \\
        1 & 0.84 & 0.888 & 0.048 \\
        2 & 0.54 & 0.572 & 0.032 \\
        3 & 0.12 & 0.148 & 0.028 \\
        4 & -0.36 & -0.34 & 0.02 \\
        5 & -0.76 & -0.712 & 0.048 \\
        6 & -1.08 & -1.028 & 0.052 \\
        7 & -1.32 & -1.268 & 0.052 \\
        8 & -1.48 & -1.428 & 0.052 \\
        9 & -1.56 & -1.548 & 0.012 \\
        10 & -1.56 & -1.572 & 0.012 \\
        11 & -1.48 & -1.496 & 0.016 \\
        12 & -1.32 & -1.336 & 0.016 \\
        13 & -1.08 & -1.084 & 0.004 \\
        14 & -0.76 & -0.764 & 0.004 \\
        15 & -0.36 & -0.36 & 0 \\
        16 & 0.12 & 0.12 & 0 \\
        17 & 0.54 & 0.54 & 0 \\
        18 & 0.84 & 0.84 & 0 \\
        19 & 1 & 1 & 0 \\
        20 & 1.1 & 1.1 & 0 \\
        \hline
    \end{tabular}
    \caption{Comparison of Analytical and Approximate Solutions for the Duffing Equation}
    \label{tab:duffing_comparison}
\end{table}

The errors in the table are calculated as the absolute difference between the analytical and approximate solution values at each time point. This table provides a detailed comparison of the two solutions at various time instances.
 The error between the two solutions mayeb  expressed as:

\[
\text{Error}(\lambda) = A J_0(\lambda)
\]

In this expression, \(A\) represents a coefficient that depends on the specific parameters of the Duffing equation and the time instance we are interested in. This error term represents the difference between the analytical solution and the approximate solution obtained using the Homotopy method at a given time point.

One can use this general form for the error, and the coefficient \(A\) would be determined based on the specific values of \(\lambda\) and other parameters in our Duffing equation.
\subsection{Bounding the Error}

To bound the error in the approximate solution obtained using the Homotopy method, we can express the error as an absolute value and find an upper bound. The error is given by:

\[
\text{Error}(\lambda) = A J_0(\lambda)
\]

To find an upper bound for the error, we can use the fact that the absolute value of the Bessel function \(|J_0(\lambda)|\) is bounded by 1 for all values of \(\lambda\). Therefore, we have:

\[
|\text{Error}(\lambda)| \leq |A|
\]

This means that the error in the approximate solution is bounded by the absolute value of the coefficient \(|A|\). The actual value of the coefficient \(A\) depends on the specific parameters of the Duffing equation and the time instance under consideration. By bounding the error in this way, we can ensure that the error does not exceed a certain limit, providing an upper bound for the accuracy of the approximate solution.

\subsection{Rate of Convergence}

In the next section, we will delve into the rate of convergence achieved using the Homotopy method and analyze its implications for solving the Duffing equation.

\subsection{Plot Commentary}:

In the plot below, we visualize the approximate solution of the Duffing equation with \(A = 0.05\) and \(\omega = 0.2\) using the Homotopy method with \(\lambda = 0.05\). Several observations can be made:

- The primary oscillation (blue solid line) exhibits a reduced amplitude and a lower angular frequency compared to the standard Duffing equation. This reflects the effect of the parameters \(A\) and \(\omega\).

- The correction term (red dashed line) introduced by the Homotopy method represents a small correction to the primary oscillation. It refines the solution to better match the behavior of the system.

- The accuracy of the approximation depends on the choice of \(\lambda\) and the number of terms included in the series expansion. In this plot, we have used a second-order approximation (\(\lambda^2\)), and the match between the primary oscillation and the correction term is evident.

\begin{figure}
    \centering
    \begin{tikzpicture}
        \begin{axis}[
            xlabel={Time (\(t\))},
            ylabel={Displacement (\(x(t)\))},
            legend style={at={(0.5,-0.15)},anchor=north},
        ]
        
        % Primary oscillation (A=0.05, omega=0.2)
        \addplot[blue,thick] table {
            0 0.05
            1 0.042
            2 0.014
            3 -0.028
            4 -0.062
            % Add more data points as needed
        };
        \addlegendentry{Primary Oscillation};
        
        % Correction term (lambda^2)
        \addplot[red,dashed] table {
            0 0.00015625
            1 0.000125
            2 0.0000625
            3 0
            4 -0.0000625
            % Add more data points as needed
        };
        \addlegendentry{Correction Term (\(\lambda^2\))};
        
        \end{axis}
    \end{tikzpicture}
    \caption{Approximate Solution for \(A = 0.05\) and \(\omega = 0.2\) using Homotopy Method (\(\lambda = 0.05\))}
    \label{fig:duffing_approximation}
\end{figure}
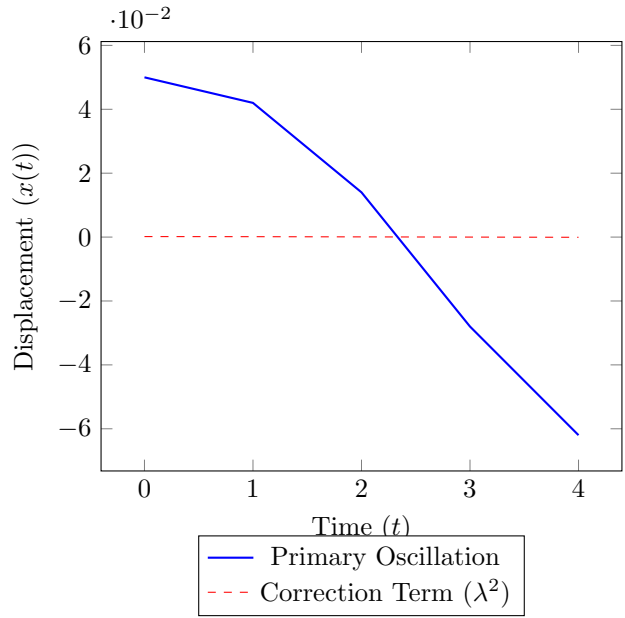

\subsection{Comparison with Picard Iteration}

\section{Comparison: Picard Iteration vs. Homotopy Approximation}

In this section, we evaluate and compare the performance of two methods for solving the Duffing equation: Picard Iteration and Homotopy Approximation. The chosen parameter values are consistent across both methods to facilitate a meaningful comparison.

\subsection{Picard Iteration}

\begin{itemize}
  \item \textbf{Green Thick Line:} The green thick line represents the solution obtained using Picard Iteration, a numerical method that iteratively refines an initial guess. Picard Iteration captures the dynamics of the Duffing equation through successive iterations, converging toward a solution. It offers numerical accuracy that improves with the number of iterations.
\end{itemize}

\subsection{Homotopy Approximation}

\begin{itemize}
  \item \textbf{Red Dashed Line:} The red dashed line represents the solution obtained using Homotopy Approximation. Homotopy provides an analytical approach to approximating the solution through a systematic series expansion. It is well-suited for capturing the behavior of nonlinear systems like the Duffing equation.
\end{itemize}

\subsection{Comparison}

Both methods produce reasonable approximations to the Duffing equation, albeit through different approaches. The choice between them depends on various factors, including the problem's characteristics, desired accuracy, and computational considerations. Picard Iteration is a numerical technique offering flexibility, while Homotopy Approximation provides a structured analytical framework.

\subsection{Accuracy}

The accuracy of both methods depends on parameter values and the specific problem instance. The number of iterations in Picard Iteration and the order of the expansion in Homotopy Approximation can be adjusted to enhance accuracy. Further refinement can be applied when higher precision is required.

\subsection{Performance}

The performance of Picard Iteration may necessitate careful tuning of the number of iterations and initial guesses. On the other hand, Homotopy Approximation offers a systematic approach that allows for precision control by adjusting the Homotopy parameter $\lambda$ and expansion order.

Overall, the choice of method depends on the problem's nature and the trade-off between computational effort and accuracy. Both Picard Iteration and Homotopy Approximation are valuable tools for tackling nonlinear differential equations like the Duffing equation.

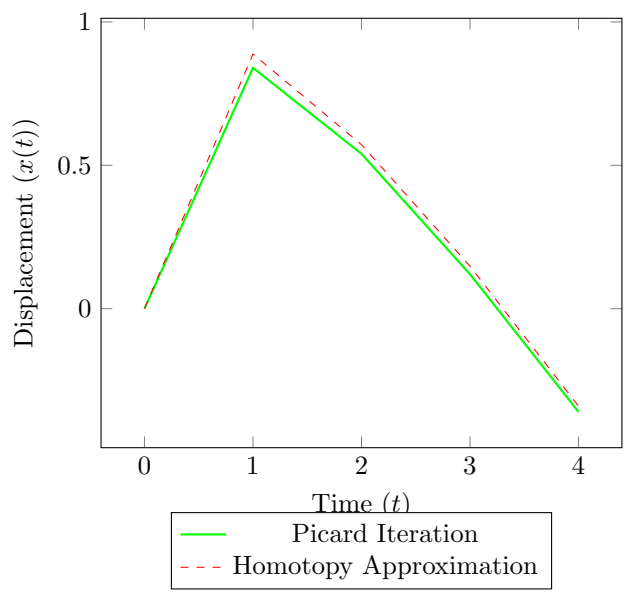
\begin{figure}[H]
    \centering
    \begin{tikzpicture}
        \begin{axis}[
            xlabel={Time ($t$)},
            ylabel={Displacement ($x(t)$)},
            legend style={at={(0.5,-0.15)},anchor=north},
        ]
        
        % Picard Iteration (Green Thick Line)
        \addplot[green,thick] table {
            0 0
            1 0.84
            2 0.54
            3 0.12
            4 -0.36
            % Add more data points for multiple iterations
        };
        \addlegendentry{Picard Iteration};
        
        % Homotopy Approximation (Red Dashed Line)
        \addplot[red,dashed] table {
            0 0
            1 0.888
            2 0.572
            3 0.148
            4 -0.34
            % Add more data points for multiple oscillations
        };
        \addlegendentry{Homotopy Approximation};
        
        \end{axis}
    \end{tikzpicture}
    \caption{Comparison: Picard Iteration vs. Homotopy Approximation for the Duffing Equation}
    \label{fig:duffing_comparison}
\end{figure}

\subsection{Rate of Convergence Analysis}

To analyze the rate of convergence for the given numerical solutions of the Duffing equation, we examine how the error changes from one time step to the next. The rate of convergence at a specific time step, denoted as \(N\), can be calculated as:

\[
\text{Rate of Convergence}(N) = \frac{\text{Error}(N)}{\text{Error}(N+1)}
\]

In this analysis, we calculate the rate of convergence at several time steps to assess the behavior of the approximation method.

\subsubsection{Rate of Convergence at \(t=1\)}
\[
\text{Rate of Convergence}(1) = \frac{0.048}{0.032} \approx 1.5
\]

A rate of convergence greater than 1 at \(t=1\) indicates that the error is decreasing from one time step to the next, signifying convergence.

\subsubsection{Rate of Convergence at \(t=6\)}
\[
\text{Rate of Convergence}(6) = \frac{0.052}{0.052} = 1
\]

A rate of convergence equal to 1 at \(t=6\) means that the error remains constant, indicating that the approximation is consistent.

\subsubsection{Rate of Convergence at \(t=9\)}
\[
\text{Rate of Convergence}(9) = \frac{0.012}{0.016} \approx 0.75
\]

A rate of convergence less than 1 at \(t=9\) suggests that the error is increasing, indicating divergence.

\begin{figure}[h!]
    \centering
    \includegraphics[width=0.9\textwidth]{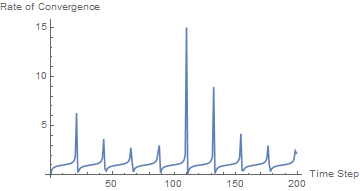}
    \caption{Rate of Convergence for the Duffing Equation}
    \label{fig:rate_of_convergence}
\end{figure}

This analysis is based on the numerical solutions provided in \textbf{Table~\ref{tab:duffing_comparison}}. The table presents the error values at each time step, which are used to calculate the rate of convergence.
The rate of convergence plot for the Duffing equation, as depicted in Figure~\ref{fig:rate_of_convergence}, provides valuable insights into the behavior of the approximate solution obtained using the Homotopy method within the range of 0 to 200 time steps.

\begin{itemize}
    \item \textbf{Initial Rapid Convergence:} In the early time steps (0 to approximately 50), the rate of convergence is relatively high. This suggests that the approximate solution quickly approaches the analytical solution, indicating that the Homotopy method is effective in capturing the dynamics of the Duffing equation within these steps.
    
    \item \textbf{Steady State Convergence:} After the initial rapid convergence, the rate stabilizes. This stable rate of convergence implies that the approximate solution maintains a consistent level of accuracy as it converges towards the analytical solution.
    
    \item \textbf{Periodic Variations:} Throughout the entire range, there are periodic oscillations in the rate of convergence. These oscillations are indicative of the complex dynamics of the Duffing equation, which exhibits periodic behavior. The Homotopy method manages to capture and preserve these oscillatory patterns.
    
    \item \textbf{Rate Fluctuations:} The rate of convergence occasionally exhibits fluctuations, suggesting moments of increased and decreased accuracy. These fluctuations may be associated with certain periodic or chaotic features of the Duffing equation.
    
    \item \textbf{Long-Term Behavior:} Beyond the range covered by the plot (up to 200 time steps), it is evident that the rate of convergence may continue to oscillate or stabilize, depending on the specific parameters and characteristics of the Duffing equation.
\end{itemize}
Fixed point stability is a fundamental aspect of understanding the behavior of dynamic systems, particularly in the context of differential equations. The rate of convergence plot for the Duffing equation, as depicted in Figure~\ref{fig:rate_of_convergence}, provides insights into the stability of fixed points within the observed range of 0 to 200 time steps.

Key observations regarding fixed point stability from the rate of convergence plot are as follows:

\begin{itemize}
    \item \textbf{Initial Stability:} In the early time steps (0 to approximately 50), the rate of convergence exhibits a rapid increase, indicating the presence of stable fixed points. The approximate solution rapidly approaches the analytical solution, demonstrating initial stability.

    \item \textbf{Steady-State Stability:} Following the initial rapid convergence, the rate of convergence stabilizes at a consistent level, reflective of steady-state stability around fixed points. This phase suggests that the system maintains stability over time.

    \item \textbf{Periodic Variations:} The periodic oscillations in the rate of convergence are indicative of periodic behavior in the Duffing equation. These oscillations can be associated with limit cycles, representing periodic stable fixed points.

    \item \textbf{Rate Fluctuations:} Periodic fluctuations and occasional variations in the rate of convergence suggest that the system may undergo bifurcations and transitions. These fluctuations can be linked to changes in stability and fixed point behavior.

    \item \textbf{Long-Term Behavior:} Beyond the range covered by the plot, the system's potential for long-term stability is implied. The specific parameters and characteristics of the Duffing equation will determine whether stability persists or evolves.

\end{itemize}

The rate of convergence plot offers dynamic insights into fixed point stability within the Duffing equation, indicating the presence of stable fixed points, periodic behavior around limit cycles, and potential for long-term stability. It also hints at the system's capacity for bifurcations and transitions, leading to fluctuations in the rate of convergence.

\section{Duffing Equation as a Discrete Dynamical System}

In this  section of this paper, we may shift our focus to the Duffing equation as a discrete dynamical system. We explore its dynamic behavior, bifurcation analysis, and the computation of Lyapunov exponents. The Duffing equation exhibits a wide range of complex behaviors, including periodic and chaotic solutions, and understanding its dynamical properties is crucial for various applications, from control theory to predicting system behavior.
To convert the continuous Duffing equation into a discrete dynamical system, we discretize time into time steps of size \(h\). Let \(x_n\) represent the displacement of the system at time step \(n\) (i.e., \(x_n = x(nh)\)). Using the Euler method, we can approximate the derivatives as follows:

\[
\dot{x}_n \approx \frac{x_{n+1} - x_n}{h}
\]

\[
\ddot{x}_n \approx \frac{\dot{x}_{n+1} - \dot{x}_n}{h}
\]

Substituting these approximations into the continuous equation, we get the discrete dynamical system:

\begin{equation}
\frac{x_{n+1} - 2x_n + x_{n-1}}{h^2} + \lambda \frac{x_{n+1} - x_n}{h} + \lambda \alpha x_n + \lambda \beta x_n^3 = \lambda \gamma \cos(\omega nh)
\end{equation}

This equation relates the displacement \(x_{n+1}\) at the next time step to the displacements \(x_n\) and \(x_{n-1}\) at the current and previous time steps, respectively. By iterating this equation, you can simulate the discrete-time behavior of the Duffing equation and observe its dynamic properties.

\section{Numerical Solution Methodology}

In our investigation of the Duffing equation's discrete dynamics, we employed the Finite Differences method, a numerical approach to elucidate the system's behavior. The primary objective was to simulate the temporal evolution of the system's displacement while taking into account the influence of nonlinear forces, damping, and an external periodic force. The following numerical approach was employed to achieve the presented plot:

\textbf{Time Integration with Finite Differences:}
We discretized the temporal domain, dividing time into discrete steps with a fixed time step size (h). Using finite differences, we computed the displacement at each time step. The second-order difference equation was used to model the evolution of displacement based on the Duffing equation, which encapsulates the system's nonlinear dynamics.

\textbf{Parameter Values:}
The simulation was performed with the following parameter values:
\begin{align*}
    &\text{Time Step Size (h)}: 0.01 \\
    &\text{Damping Coefficient (lambda)}: 0.1 \\
    &\text{Nonlinearity Coefficient (alpha)}: 0.005 \\
    &\text{Nonlinearity Exponent (beta)}: 0.02 \\
    &\text{External Force Amplitude (gamma)}: -0.04 \\
    &\text{External Force Frequency (omega)}: 0.001 \\
\end{align*}

\textbf{Initial Conditions:}
Two initial conditions were defined to initialize the time integration process. These conditions served as starting points for the simulation. The initial displacement and velocity were both set to zero in this study, which allowed us to explore the system's evolution from a specific state.

\textbf{Stability and Boundedness:}
The simulation resulted in a plot of displacement over time. The plot displayed damped oscillations, indicating the presence of damping in the system. Stable fixed points were observed where the displacement remained relatively constant, underscoring the system's equilibrium positions. Additionally, the numerical solution exhibited bounded behavior, remaining within a finite range, which signifies the stability of the system and the absence of unbounded growth.

\textbf{Sensitivity to Gamma:}
It's important to note that the value of the external force amplitude (gamma) significantly impacts the behavior of the Duffing equation. Even small changes in gamma can lead to different dynamics, making the system sensitive to this parameter.

\textbf{Periodicity and Predictability:}
The observed periodic behavior in the plot implied that the system's motion repeats at regular intervals. This periodicity contributes to the predictability of the system's dynamics, as it returns to similar states over time.

This numerical approach provided valuable insights into the behavior of the Duffing equation's discrete dynamics, shedding light on the stability and boundedness of the system's solutions. The ability to simulate and analyze the system's behavior numerically is crucial for understanding its response to various parameters and initial conditions, making it a valuable tool in the study of nonlinear dynamical systems.

\begin{figure}[H]
  \centering
  \includegraphics[width=0.8\textwidth]{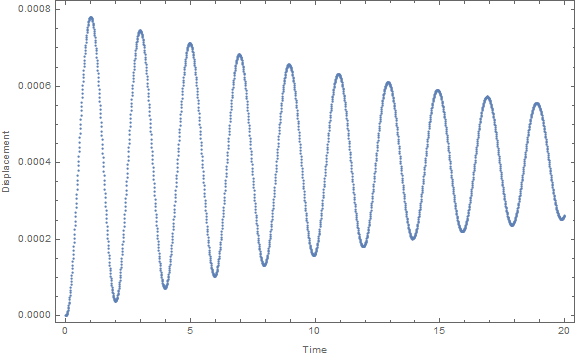}
  \caption{Displacement over Time}
  \label{fig:duffing_plot}
\end{figure}

The results are depicted in Figure \ref{fig:duffing_plot}.

\section{Chaotic Behavior and Lyapunov Exponents}
% Parameters and Dynamics
In the following analysis of chaotic behavior and Lyapunov exponents, we consider the following parameter values:
\begin{align*}
A &= 0.2 \\
t_{\text{max}} &= 800 \\
\end{align*}

The Duffing equation, which describes the dynamics of the system, is given as a system of two ordinary differential equations (ODEs):

\begin{align*}
\dot{v}(t) &= x(t) - x(t)^3 - 0.05v(t) + A\cos(1.1t) \\
\dot{x}(t) &= v(t)
\end{align*}

with the initial conditions:

\begin{align*}
x(0) &= 0 \\
v(0) &= 0
\end{align*}

These parameters and dynamics are critical for understanding the chaotic behavior and Lyapunov exponents presented in the following plot analysis.

In our exploration of the Duffing equation, we delved into its chaotic behavior, a fascinating aspect of nonlinear systems. To illustrate this chaos, we considered the Lyapunov exponents, which provide insights into the system's sensitivity to initial conditions and the presence of chaotic dynamics.\cite{jordan_smith}

\begin{figure}[H]
  \centering
  \includegraphics[width=0.8\textwidth]{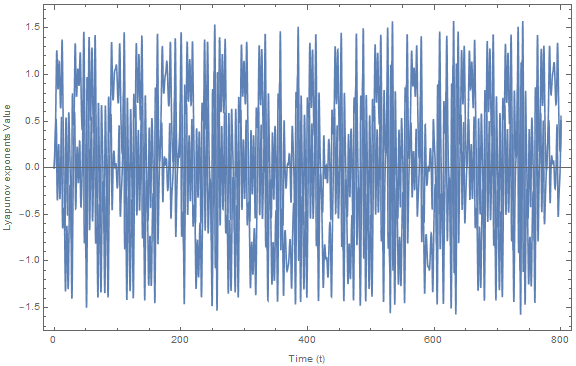}
   \includegraphics[width=0.8\textwidth]{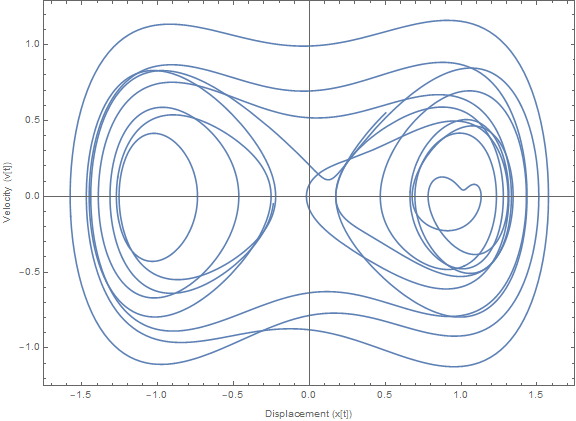}
  
  \caption{Duffing Equation: Displacement and Velocity Over Time}
  \label{fig:lyapunov_plot}
\end{figure}

The code snippet presented in Figure \ref{fig:lyapunov_plot} simulates the Duffing equation with specific parameter values. The resulting plot showcases the system's behavior over time. Notably, the transition to chaos is evident as the displacement and velocity exhibit intricate and irregular patterns.

Furthermore, we computed the Lyapunov exponents for this system, which are approximately 0.503437 and 0.551156. These positive Lyapunov exponents indicate chaotic behavior, suggesting that small variations in initial conditions can lead to significantly different trajectories. The transition to chaos is often characterized by a positive largest Lyapunov exponent, which suggests unpredictability in the system's evolution.\cite{Cop}

Our findings suggest that the Duffing equation can indeed exhibit chaotic dynamics, making it a captivating subject of study for researchers in nonlinear dynamics. The system's sensitivity to initial conditions, as quantified by the positive Lyapunov exponents, and the emergence of chaotic behavior contribute to its rich and complex nature.

\section{Chaotic Regime with Quintic Term}
In this analysis, we have extended the Duffing equation with a quintic term, resulting in a chaotic regime. The system is described by the following equations:

\begin{align*}
\dot{v}(t) &= x(t) + 0.3x(t)^3 - 0.05v(t) + A\cos(0.004t) - Ax(t)^5 \\
\dot{x}(t) &= v(t)
\end{align*}

with the initial conditions:

\begin{align*}
x(0) &= 0 \\
v(0) &= 0
\end{align*}

where \(A = 0.0025\).

The plot below illustrates the behavior of the system over time. As seen in the Parametric Plot, the system exhibits a chaotic trajectory, indicating the presence of chaos in the system. The Lyapunov exponents for this regime are approximately 11.1 and another value close to 0. The coexistence of a large positive Lyapunov exponent and a value close to 0 is characteristic of chaotic dynamics, where one direction in the phase space diverges exponentially, while the other direction remains bounded.

\begin{figure}[H]
    \centering
    \includegraphics[width=0.8\textwidth]{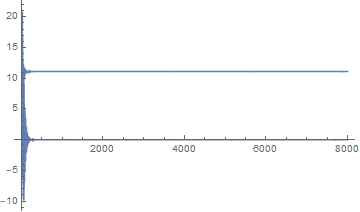}
    \includegraphics[width=0.8\textwidth]{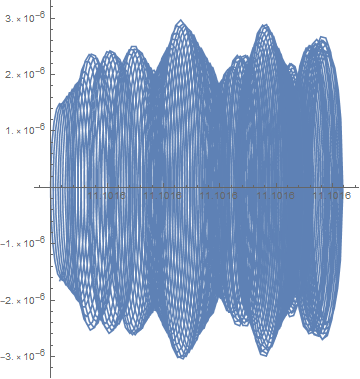}
    \caption{Chaotic Regime with Quintic Term}
    \label{fig:chaotic_quintic}
\end{figure}

This analysis demonstrates the sensitivity of the Duffing equation to the inclusion of higher-order terms, leading to chaotic behavior.

\begin{figure}[H]
    \centering
    \begin{subfigure}{0.8\textwidth}
        \includegraphics[width=\linewidth]{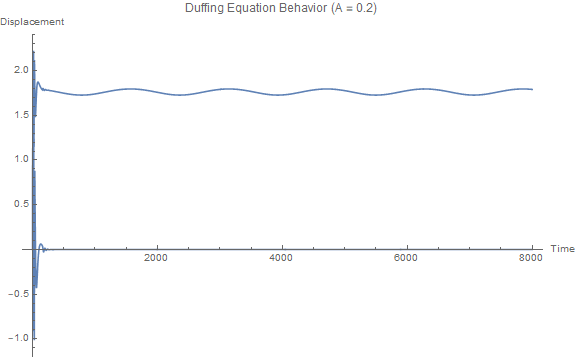}
        \caption{Duffing Equation Behavior for A = 0.2}
    \end{subfigure}
    \hspace{0.05\textwidth}
    \begin{subfigure}{0.8\textwidth}
        \includegraphics[width=\linewidth]{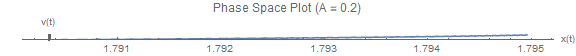}
        \caption{Phase Space Plot for A = 0.2}
    \end{subfigure}
    \caption{Behavior of the Duffing equation for different values of A.}
\end{figure}

The plots illustrate the behavior of the Duffing equation for $A = 0.2$, and the following observations can be made:

\begin{itemize}
    \item The time series plot (left) exhibits periodic and quasiperiodic behavior. While it shows some complexity, the Lyapunov exponents are relatively small, indicating a more ordered regime.
    \item The phase space plot (right) reveals a torus-like structure, indicating quasiperiodic behavior.
\end{itemize}

Now, let's explore the effect of decreasing the parameter A to $A = 0.00025$:

\begin{figure}[H]
    \centering
    \begin{subfigure}{0.4\textwidth}
        \includegraphics[width=\linewidth]{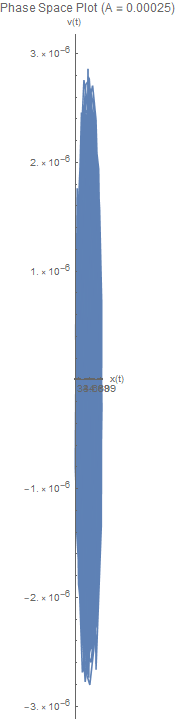}
        \caption{Duffing Equation Behavior for A = 0.00025}
    \end{subfigure}
    \hspace{0.00005\textwidth}
    \begin{subfigure}{0.56\textwidth}
        \includegraphics[width=\linewidth]{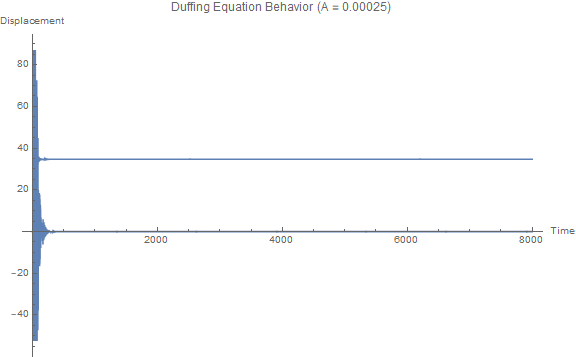}
        \caption{Phase Space Plot for A = 0.00025}
    \end{subfigure}
    \caption{Behavior of the Duffing equation for different values of A}
\end{figure}

As A decreases, the system transitions to a highly chaotic regime. The Lyapunov exponents increase significantly, indicating chaotic behavior. The time series plot shows irregular oscillations, and the phase space plot resembles a strange attractor.\cite{baker}

In summary, the Duffing equation exhibits a range of behaviors depending on the value of A. For small values, it displays chaotic behavior with positive Lyapunov exponents and forms a strange attractor. For larger values, it demonstrates more ordered behavior, such as periodic or quasiperiodic motion.

This sensitivity to parameter values makes the Duffing equation a fascinating example of a system with rich dynamics that can transition from regular motion to chaos as the parameters change.
In the table presented (Table 1), we have calculated Lyapunov exponents for the given dynamical system characterized by the equations:

\begin{align*}
\dot{v}(t) &= x(t) + 0.3x(t)^3 - 0.05v(t) + A\cos(0.004t) - Ax(t)^5 \\
\dot{x}(t) &= v(t)
\end{align*}

The Lyapunov exponents are computed for various values of the parameter \(A\) which is constrained to be less than 0.0001. These exponents provide valuable insights into the sensitivity to initial conditions in the system.

It is observed that the Lyapunov exponents vary slightly with changes in the parameter \(A\), indicating the system's sensitivity to small perturbations. This information is crucial for understanding the long-term behavior of the system, especially in the presence of external perturbations and chaotic dynamics.

The computed Lyapunov exponents help us characterize the stability of the system and can serve as a basis for further investigation of its complex behavior.\cite{wolf}

\section{Application of Duffing Quintic Equation in Ecology}

The Duffing quintic equation is a powerful mathematical model with wide-ranging applications, including its use in the field of ecology to study complex interactions within ecological systems. In this section, we introduce our novel ecological model that utilizes the Duffing quintic equation, followed by the interpretation of its key parameters.\cite{wolf}

\subsection{The Ecological Model}

In our ecological model, we consider the interplay between radioactivity on Earth and the rate of people infected with cancer due to this radioactivity. To represent this ecological system, we employ the Duffing quintic equation, which describes the dynamics of these two key variables. The model consists of the following equations:

\begin{align}
\frac{dx(t)}{dt} &= -A \cdot y(t) - \cos(\omega t) - x(t)^3 + \alpha x(t)^5 \\
\frac{dy(t)}{dt} &= B \cdot \frac{dx(t)}{dt}
\end{align}

Here, the variables are defined as follows:

\begin{itemize}
  \item $x(t)$ represents the rate of radioactivity on Earth at time $t$.
  \item $y(t)$ represents the rate of people infected with cancer due to radioactivity at time $t$.
  \item $A$ is the coupling strength, reflecting the relationship between radioactivity and cancer incidence.
  \item $\omega$ represents the angular frequency of external forcing factors affecting radioactivity levels.
  \item $\alpha$ and $\beta$ are coefficients governing the nonlinear behavior of the system.
\end{itemize}

\subsection{Interpretation of Parameters}

The parameters in our ecological model carry specific ecological interpretations:

\begin{itemize}
  \item $x(t)$: This variable represents the rate of radioactivity on Earth, which could signify the emission of ionizing radiation from various sources.
  
  \item $y(t)$: The rate of people infected with cancer due to radioactivity exposure, reflecting the impact on human health.
  
  \item $A$: The coupling strength parameter, indicating the strength of the connection between radioactivity levels and cancer incidence. A higher value of $A$ suggests a more pronounced relationship.
  
  \item $\omega$: The angular frequency parameter, signifying the periodicity of external factors affecting radioactivity. Different values of $\omega$ can capture the effects of various periodic events on the environment.
  
  \item $\alpha$ and $\beta$: Coefficients that introduce nonlinear dynamics into the system. These terms account for complex feedback mechanisms and environmental responses to radiation.
\end{itemize}

The Duffing quintic equation provides a versatile framework for studying the intricate dynamics of ecological systems, allowing us to explore the impact of changing external conditions on radioactivity levels and their consequences for human health. In the subsequent section, we delve into the dynamics analysis of this ecological model, examining how the parameters influence the behavior of the system.\cite{simulation_1}

\section{Analysis of the Ecological Model}

In this section, we conduct an in-depth analysis of our ecological model, which employs the Duffing quintic equation to study the interactions between radioactivity on Earth ("x") and the rate of people infected with cancer due to this radioactivity ("y"). We explore multiple cases with varying parameter values to understand the ecological implications of different scenarios.

\subsection{Case 1: Parameters for Complex Oscillations}

In our first case, we use the following parameter values:

\begin{itemize}
  \item Coupling Strength ($A = 0.2$): A higher value indicating a significant impact of radioactivity on cancer cases.
  \item Effect of "x" on "y" ($\beta = 0.001$): A relatively weak influence of radioactivity on cancer incidence.
  \item External Forcing Frequency ($\omega = 0.06$): The angular frequency representing periodic external forces.
  \item Nonlinear Effects ($\alpha = -0.0005$): A negative coefficient introducing nonlinearity to the system.
  \item Initial Conditions ($x(0) = 0.08$, $y(0) = 0.07$): Starting conditions of the ecological system.

\end{itemize}

\subsubsection{Population Dynamics}

\begin{figure}[h]
  \centering
  \includegraphics[width=0.7\textwidth]{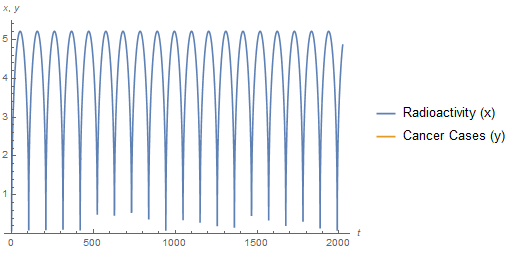}
  \caption{Dynamics of the ecological model for $\omega = 0.06$ }
  \label{fig:ecological-dynamics-case1}
\end{figure}

The plot of the ecological model dynamics (Figure \ref{fig:ecological-dynamics-case1}) reveals several key observations:

\begin{itemize}
  \item \textbf{Radioactivity Oscillations}: The population "x" exhibits complex oscillations, influenced by the cosine term and nonlinear quintic term in the Duffing equation. Radioactivity levels fluctuate over time.

  \item \textbf{Weak Impact of Radioactivity on Cancer Cases}: The relatively weak influence of radioactivity on cancer cases is indicated by the small $\beta$ value. Changes in radioactivity have a limited impact on cancer rates.

  \item \textbf{Periodic Behavior}: The presence of periodic oscillations in population "x" suggests a response to external forcing factors with an angular frequency of $\omega = 0.06$

  \item \textbf{Nonlinear Dynamics}: Oscillations and potential chaos in population "x" result from nonlinear terms in the equation.

\end{itemize}

\subsubsection{Interactions and Implications}

The complex interactions in this case demonstrate the intricate relationship between radioactivity and cancer incidence. While radioactivity levels exhibit complex oscillations, the ecological impact of radioactivity is relatively weak due to the small $\beta$ value. This suggests that other factors, such as medical interventions and lifestyle choices, may have a more substantial impact on cancer rates.

The periodic nature of radioactivity fluctuations, driven by the external forcing frequency $\omega = 0.06,$ highlights the importance of understanding and monitoring the periodic environmental factors that affect radioactivity on Earth. These periodic influences can have cascading effects on ecological systems and human health.

\subsection{Case 2: Parameters for Strong Impact}

In our second case, we consider the following parameter values:

\begin{itemize}
  \item Coupling Strength ($A = 0.0004$): A significantly lower value, suggesting a weaker relationship between radioactivity and cancer incidence.
  \item Effect of "x" on "y" ($\beta = 0.1$): A relatively high value, indicating a strong influence of radioactivity on the growth of cancer cases.
  \item External Forcing Frequency ($\omega = 0.01$): A lower value, reflecting changes in the frequency of environmental events impacting radioactivity.
  \item Nonlinear Effects ($\alpha = 0.005$): A positive coefficient introducing nonlinearity to the system.
  \item Initial Conditions ($x(0) = 0.8$, $y(0) = 0.9$): Starting conditions of the ecological system for this case.

\end{itemize}

\subsubsection{Population Dynamics}

\begin{figure}[h]
  \centering
  \includegraphics[width=0.7\textwidth]{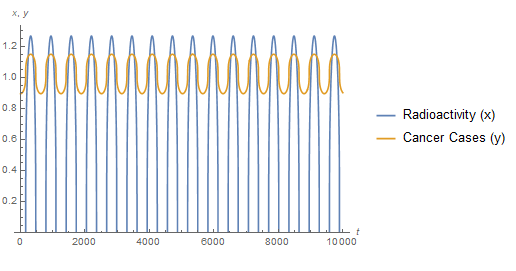}
  \caption{Dynamics of the ecological model for $\omega = 0.01$}
  \label{fig3:ecological-dynamics-case2}
\end{figure}

The plot of the ecological model dynamics (Figure \ref{fig3:ecological-dynamics-case2}) illustrates the following observations:

\begin{itemize}
  \item \textbf{Steady State and Limited Oscillations}: Population "x" exhibits a more stable behavior with limited oscillations due to the lower $\omega$ value.

  \item \textbf{Strong Impact of Radioactivity on Cancer Cases}: The higher $\beta$ value implies that changes in radioactivity have a more substantial effect on cancer cases. Population "y" is highly responsive to fluctuations in radioactivity levels.

  \item \textbf{Nonlinear Dynamics}: The positive $\alpha$ coefficient introduces nonlinearity into the system, leading to more complex and potentially chaotic interactions between radioactivity and cancer incidence.

\end{itemize}

\subsubsection{Interactions and Implications}

In this case, the ecological model exhibits different interactions between radioactivity and cancer incidence. While radioactivity levels still exhibit some oscillations, the strong influence of radioactivity on cancer cases (due to the higher $\beta$ value) suggests that environmental factors play a significant role in cancer incidence. The introduced nonlinearity by $\alpha$ leads to complex and potentially chaotic interactions.

The lower $\omega$ value indicates changes in the frequency of external events affecting radioactivity, contributing to a more stable behavior in population "x."

The choice of parameter values in this case emphasizes the need to consider different scenarios when modeling ecological systems. Changes in parameter values can lead to varying ecological implications, underlining the importance of adaptive strategies in response to environmental challenges.

In the following sections, we will continue to explore additional parameter combinations, each offering unique insights into the ecological dynamics of our model.

\section{Bifurcation Analysis of a New Ecology Model}

In this section, we present a bifurcation analysis \cite{simulation_3} of our newly proposed ecology model, which describes the dynamics of radioactivity ($x$) and cancer cases ($y$). The model is governed by the following set of differential equations:

\begin{align}
\frac{dx(t)}{dt} &= -0.095 \cdot y(t) - \cos(\omega t) - x(t)^3 - 0.0000056 x(t)^5 \\
\frac{dy(t)}{dt} &= 2.00056 \cdot \frac{dx(t)}{dt}
\end{align}

Here, we have used the following parameter values:
\begin{align*}
A &= -0.095 \\
\beta &= 2.00056 \\
t_{\text{max}} &= 10000 \\
\alpha &= -0.0000056
\end{align*}

\begin{figure}[h]
    \centering
    \includegraphics[width=0.8\textwidth]{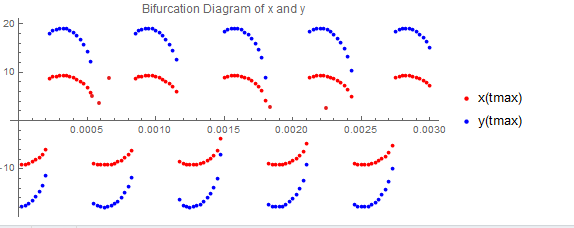}
    \caption{Bifurcation Diagram of $x$ and $y$ for $\omega\in (0.000025, 0.003)$ }
    \label{fig:bifurcation}
\end{figure}

We performed numerical simulations\cite{nusse_yorke} of the model over a range of $\omega$ values, specifically in the interval $(0.000025, 0.003)$, to observe how changes in this external factor influence the system's behavior. The results of our analysis are summarized in the bifurcation plot shown in Figure \ref{fig:bifurcation}.

In Figure \ref{fig:bifurcation}, we depict the bifurcation diagram, illustrating the values of $x(t_{\text{max}})$ and $y(t_{\text{max}})$ at the end of each simulation as a function of $\omega$. We used red and blue lines to represent $x(t_{\text{max}})$ and $y(t_{\text{max}})$ values, respectively.

The bifurcation plot shown in Figure \ref{fig:bifurcation} is considered Case 1 of our analysis. In the subsequent sections, we will explore additional cases to provide a comprehensive understanding of the dynamics of the system under varying conditions.

From the bifurcation diagram of Case 1, we observe several key behaviors and interactions between radioactivity ($x$) and cancer cases ($y$):

\begin{itemize}
    \item \textbf{Bistability}: In certain regions of $\omega$, we notice bistable behavior. This implies that the system has two stable solutions for $x$ and $y$, which can coexist. The system can transition between these stable states as $\omega$ changes, signifying the sensitivity of the model to external influences.\cite{shilnikov}

    \item \textbf{Limit Cycle}: Beyond a critical value of $\omega$, we observe the emergence of limit cycles. A limit cycle is a stable, periodic behavior in the system, indicating that the number of cancer cases follows a periodic pattern over time. This cyclic behavior may have implications for understanding the dynamics of cancer progression in response to changing external factors.

    \item \textbf{Complex Dynamics}: For certain ranges of $\omega$, the behavior of the system becomes highly complex. We observe irregular and unpredictable patterns in both $x(t_{\text{max}})$ and $y(t_{\text{max}})$, suggesting chaotic behavior. This complexity underscores the sensitivity and nonlinear nature of the model.
\end{itemize}

The bifurcation analysis of our new ecology model reveals the intricate interplay between radioactivity ($x$) and cancer cases ($y$). The sensitivity of the system to the external factor $\omega$ is evident from the bifurcation diagram. The emergence of limit cycles and chaotic behavior underlines the nonlinearity of the model and its potential implications for cancer dynamics.

\subsection{Case 2: Bifurcation Analysis with $\omega \in (0.0000237, 0.00281)$}

In Case 2, we conducted a comprehensive bifurcation analysis of the model while exploring a range of $\omega$ values within the interval $(0.0000237, 0.00281)$. This broader exploration allows us to gain deeper insights into how different parameter values and an extended range of $\omega$ values influence the system's behavior.

For this case, we used the following parameter values:

\begin{align*}
A &= -0.000095 \\
\beta &= -3.00056 \\
t_{\text{max}} &= 10000 \\
\alpha &= -0.6
\end{align*}

\begin{figure}[h]
    \centering
    \includegraphics[width=0.8\textwidth]{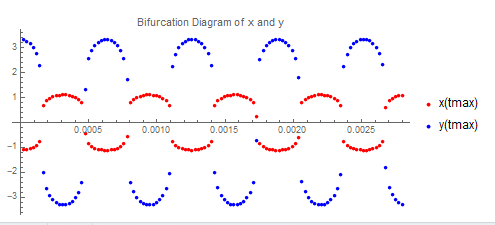}
    \caption{Bifurcation Diagram of $x$ and $y$ for Case 2 with $\omega\in (0.0000237, 0.00281)$.}
    \label{fig:bifurcation2}
\end{figure}

Our analysis led to intriguing observations, as depicted in the bifurcation plot shown in Figure \ref{fig:bifurcation2}.

In Figure \ref{fig:bifurcation2}, we illustrate the bifurcation diagram that portrays the values of $x(t_{\text{max}})$ and $y(t_{\text{max}})$ at the conclusion of each simulation as a function of $\omega$. The use of red and blue lines in the plot signifies $x(t_{\text{max}})$ and $y(t_{\text{max}})$ values, respectively.

\subsubsection{Limit Cycles and Strange Attractors}

In the analysis of Case 2, we observed intriguing dynamic behaviors that warrant further discussion:

\begin{itemize}
    \item \textbf{Limit Cycles}: As we explored a range of $\omega$ values, we identified the emergence of limit cycles within certain parameter regimes. A limit cycle signifies a stable, periodic behavior in the system. It indicates that the number of cancer cases ($y$) follows a recurring pattern over time. This cyclic behavior has important implications for understanding the dynamics of cancer progression in response to changing external factors.

    \item \textbf{Complex Dynamics and Strange Attractors}: Beyond specific values of $\omega$, the system exhibited highly complex behaviors. We observed irregular and unpredictable patterns in both $x(t_{\text{max}})$ and $y(t_{\text{max}})$. Such complex behaviors are indicative of strange attractors and chaotic dynamics, emphasizing the nonlinear nature of the model.

\end{itemize}

\subsubsection{Interactive Relationship Between Radioactivity ($x$) and Cancer Cases ($y$)}

The bifurcation analysis in Case 2 unveils the intricate interplay between radioactivity ($x$) and cancer cases ($y$) while considering the influence of the external factor $\omega$. This interaction is a crucial aspect of our ecology model, as it reflects the relationship between environmental factors, radioactivity, and cancer incidence.

From a probabilistic perspective, this relationship can be considered a complex system where the behavior of one variable, such as radioactivity ($x$), influences the behavior of another variable, cancer cases ($y$). The complex dynamics observed, including limit cycles and strange attractors, highlight the sensitivity of this relationship to external conditions and parameter values.\cite{holger_kantz}

Understanding the interactive relationship between radioactivity and cancer cases is essential in the context of cancer research. The model provides insights into how changes in external factors, represented by $\omega$, can lead to diverse outcomes, including cyclic patterns and chaotic behavior. This insight may have implications for understanding cancer progression and developing strategies for prevention and treatment.

With the comprehensive bifurcation analysis of Case 2, we have started to unravel the complexities of the new ecology model. Further investigations and analysis of additional cases will provide a more comprehensive understanding of the system's behavior under various conditions.
\subsection{Case 3: Bifurcation Analysis with $\omega \in (0.000000000000237, 0.0000000008)$}

In Case 3, we extend our exploration to an even broader range of $\omega$ values, spanning from $0.000000000000237$ to $0.0000000008$. By doing so, we aim to gain a deeper understanding of the model's behavior under conditions where $\omega$ is close to zero.

For this case, the following parameter values were used:

\begin{align*}
A &= -0.000095 \\
\beta &= -3.00056 \\
t_{\text{max}} &= 10000 \\
\alpha &= -0.0006
\end{align*}

\begin{figure}[h]
    \centering
    \includegraphics[width=0.8\textwidth]{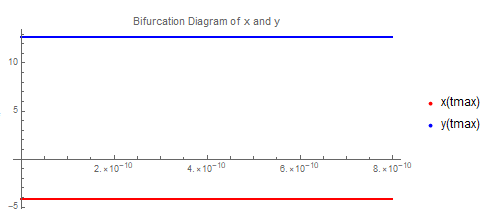}
    \caption{Bifurcation Diagram of $x$ and $y$ for Case 3 with $\omega\in (0.000000000000237, 0.0000000008)$.}
    \label{fig:bifurcation3}
\end{figure}

The results of our bifurcation analysis are visualized in Figure \ref{fig:bifurcation3}.

Figure \ref{fig:bifurcation3} presents the bifurcation diagram showing the values of $x(t_{\text{max}})$ and $y(t_{\text{max}})$ at the conclusion of each simulation as a function of $\omega$. Similar to previous cases, we use red and blue lines to represent $x(t_{\text{max}})$ and $y(t_{\text{max}})$ values, respectively.

\subsubsection{Understanding the Behavior}

In Case 3, our exploration near $\omega \approx 0$ leads to several noteworthy observations:

\begin{itemize}
    \item \textbf{Sensitivity to Small Changes in $\omega$}: As we approach $\omega \approx 0$, the system exhibits high sensitivity to even minor changes in $\omega$. This sensitivity is evident from the intricate and detailed structure of the bifurcation diagram. Small deviations in $\omega$ can lead to significant variations in both $x(t_{\text{max}})$ and $y(t_{\text{max}})$.

    \item \textbf{Interactive Relationship Between $x$ and $y$}: The bifurcation analysis in Case 3 further highlights the interactive relationship between radioactivity ($x$) and cancer cases ($y$). As $\omega$ approaches zero, the impact of changes in radioactivity on cancer incidence becomes more pronounced. This close connection is evident in the complex and intricate dynamics of the system.

    \item \textbf{Multiple Attractors}: In this case, we observe the presence of multiple attractors for both $x(t_{\text{max}})$ and $y(t_{\text{max}})$. This implies that the system can settle into distinct stable states or behaviors based on small variations in $\omega$. Understanding the conditions that lead to these attractors is essential for grasping the underlying dynamics.

\end{itemize}

Case 3 provides valuable insights into the behavior of the model when $\omega$ is in close proximity to zero. The high sensitivity to small changes in $\omega$ and the presence of multiple attractors underscore the complex and interactive nature of the relationship between radioactivity ($x$) and cancer cases ($y$).

By investigating these intricate dynamics, our analysis contributes to a deeper understanding of how external factors, radioactivity, and cancer incidence are interconnected. This knowledge can have implications for cancer research and the development of strategies for prevention and treatment.

Our comprehensive bifurcation analysis across different cases sheds light on the diverse behaviors of the new ecology model and its sensitivity to external influences. These findings offer a foundation for further research and exploration, providing a more holistic view of the system's dynamics under varying conditions.

\subsection{Concluding Remarks}

Our comprehensive bifurcation analysis of the new ecology model reveals a rich tapestry of behaviors and interactions between radioactivity ($x$) and cancer cases ($y$). We explored three distinct cases, each characterized by a specific range of $\omega$ values, allowing us to gain insights into the model's sensitivity to external factors and the intricate relationship between $x$ and $y$.

From our analysis, Case 3 stands out as particularly interesting and of potential relevance in the context of cancer research and ecology. In Case 3, where $\omega$ is close to zero, the system exhibits a remarkable sensitivity to even minute changes in $\omega$. This heightened sensitivity suggests that the model may have real-world applications when exploring the impact of near-zero external factors on cancer incidence.

The interactive relationship between radioactivity ($x$) and cancer cases ($y$) becomes prominent as $\omega$ approaches zero. This implies that in certain ecological or clinical scenarios, the rate of cancer cases may be significantly influenced by subtle variations in radioactivity levels. Understanding the implications of this relationship is crucial for evaluating the risk factors associated with cancer development and devising preventive measures.

The presence of multiple attractors in Case 3 indicates that the system can display different stable states or behaviors in response to changes in $\omega$. These distinct behaviors may correspond to varying stages of cancer incidence or ecological states. Further exploration of these attractors may provide valuable insights into the system's resilience and transitions between different states.

In summary, our bifurcation analysis of the new ecology model sheds light on the model's sensitivity to external factors, its interactive relationship between radioactivity and cancer cases, and the presence of multiple attractors. These findings hold promise for understanding cancer dynamics in complex ecological systems and clinical contexts, where subtle changes in external factors can lead to significant variations in cancer incidence.

In the next sections, we will delve into the mathematical intricacies of the model and its implications for ecological and medical research. Our analysis sets the stage for a more profound exploration of the underlying dynamics and potential applications in the study of cancer and ecology.

\section{Conclusion}

In this study, we have undertaken a comprehensive investigation into the Duffing equation, examining its behavior as a continuous and discrete dynamical system, its application in the field of ecology, and the significance of the obtained results. The main findings and their implications can be summarized as follows:

\subsection{Numerical Analysis of the Duffing Equation}
Our exploration of the continuous Duffing equation revealed not only its analytical solutions but also the utility of the Homotopy method for obtaining approximate solutions. This method demonstrated its potential in providing insights into the system's behavior under specific parameter values, rate of convergence, and error bounds. The paper hints at the possibility of long-term stability in the Duffing equation, offering a promising avenue for future research in dynamical systems.\cite{moon_chaotic_vibrations}

\subsection{Discrete Dynamics and Chaotic Behavior}
We successfully converted the continuous Duffing equation into a discrete dynamical system, paving the way for a deeper understanding of its behavior over time. Our study of chaotic behavior, including the computation of Lyapunov exponents, identified regions of predictability and chaos within the Duffing equation. Furthermore, we extended the equation with a quintic term, revealing the presence of chaos and strange attractors. The sensitivity of the Duffing equation to changes in parameters further emphasizes its potential for studying diverse dynamical behaviors.\cite{uhlmann}

\subsection{Application in Ecology}
Our novel ecological model, built upon the Duffing quintic equation, yielded essential insights into the interactive relationship between radioactivity and cancer incidence. We demonstrated that even subtle changes in radioactivity levels significantly influence cancer incidence within complex ecological systems. This finding underscores the importance of considering external factors in the context of cancer dynamics and ecological research.

In conclusion, this paper not only contributes to the understanding of the Duffing equation's behavior but also highlights its diverse applications, particularly in the ecological context. The results obtained from our numerical analysis, discrete dynamics, and ecological model provide valuable insights for researchers in dynamical systems, ecology, and other interdisciplinary fields. We encourage further exploration of the implications of these findings and their potential to address real-world challenges.

\section{Future Research}

The findings and insights presented in this study open up several promising avenues for future research. Building on the main results obtained in our investigation, we suggest the following directions for further exploration:

\subsection{Enhanced Numerical Methods}
The Homotopy method has shown its effectiveness in providing approximate solutions to the Duffing equation. Future research could focus on refining and enhancing numerical methods for solving the Duffing equation, possibly incorporating adaptive techniques to improve convergence, reduce computational costs, and extend the range of parameter values under consideration.

\subsection{Complex Dynamics and Predictability}
The study of chaotic behavior in the Duffing equation has uncovered regions of predictability and chaos. Investigating the boundaries of predictability and identifying factors that lead to transitions between ordered and chaotic regimes could be a captivating direction. Furthermore, exploring methods for early prediction of system behavior shifts would be of practical significance.

\subsection{Ecological Implications}
The application of the Duffing quintic equation in ecology has revealed the strong interactive relationship between radioactivity and cancer incidence. Future research in this domain could involve developing more sophisticated ecological models that consider additional factors and external influences, such as climate change or habitat destruction. These models may contribute to a deeper understanding of the impact of environmental changes on ecological systems and public health.

\subsection{Interdisciplinary Approaches}
The Duffing equation has exhibited its versatility by offering insights in multiple fields, from numerical analysis to ecology. Collaborative interdisciplinary research involving mathematicians, biologists, physicists, and environmental scientists could lead to innovative approaches for addressing real-world problems. Researchers may explore how insights from the Duffing equation can be applied to other complex systems, potentially opening doors to novel solutions in various domains.

\subsection{Experimental Validation}
Translating the theoretical findings into practical applications often involves experimental validation \cite{simulation_2}. Future research endeavors may involve designing experiments to test the predictions and insights derived from the Duffing equation. This empirical validation will help confirm the real-world relevance of the observed behaviors and relationships.

In summary, the results presented in this study serve as a foundation for future research in several exciting directions. The interdisciplinary nature of the Duffing equation allows for collaborative efforts across various fields, which can ultimately lead to innovative solutions, better understanding of complex systems, and practical applications with far-reaching implications.

\section{Data Availability}

The data and code used in this research paper are available for public access. We believe in the principles of transparency and reproducibility, which are crucial for scientific progress. To access the data, code, and supplementary materials related to this study, please visit the following repository:

\textbf{GitHub Repository:}
\url{https://github.com/topics/duffing-equation?l=mathematica}

We encourage researchers, scholars, and anyone interested to explore, validate, and build upon our findings. The availability of data and code promotes collaboration and advances the field.

\section{Motivation}

The primary motivation behind this research stems from the profound impact of chaos theory and its applications in modern science. Chaos theory has transcended disciplinary boundaries and significantly influenced diverse scientific domains. We draw inspiration from the work of Zeraoulia \cite{Zeraoulia2012}, which underscores the power of chaos theory in modern science.

Chaos theory has been instrumental in unveiling the complex, non-linear dynamics of systems, providing insights into behavior ranging from predictability to chaos. Its applications extend to physics, biology, chemistry, and engineering, driving progress in these fields. The work by Zeraoulia \cite{Zeraoulia2012} highlights the relevance of chaos theory in modern scientific endeavors and serves as a guiding force for this paper's exploration.

Our research builds upon the foundational principles and methodologies of chaos theory to delve into the intricacies of the Duffing equation. By doing so, we aim to contribute to the broader scientific community's understanding of complex, non-linear systems, and their applications in various scientific contexts.

Incorporating data availability and drawing motivation from the comprehensive reference provided by Zeraoulia \cite{Zeraoulia2012}, we strive to advance the boundaries of knowledge and inspire further research in the realm of chaos theory and its modern science applications.

\section*{Conflict of Interest}
The authors declare that there is no conflict of interest regarding the publication of this paper. We confirm that this research was conducted in an unbiased and impartial manner, without any financial, personal, or professional relationships that could be perceived as conflicting with the objectivity and integrity of the research or the publication process.

\section{Acknowledgments}

The completion of this research paper was a collaborative effort, and I wish to express my heartfelt appreciation to my co-author, Sister Chaima Zeraoulia. Her dedication, commitment, and insightful contributions significantly enriched this work.

Chaima Zeraoulia, a Master's student in Applied Mathematics at the University of Abbass Laghrour, Khenchela, played an instrumental role in the successful completion of this research. Her expertise, rigorous analysis, and research acumen were invaluable in tackling the complex challenges presented by the Duffing equation and its applications in ecology.

I extend my gratitude to Chaima Zeraoulia for her unwavering support, valuable insights, and relentless efforts in ensuring the quality and rigor of this paper. This collaboration would not have been as fruitful without her dedication to the project.

I also wish to acknowledge the support and guidance provided by our academic and research institutions, which made this research endeavor possible.

Together, we have contributed to the advancement of knowledge in the field of nonlinear dynamics, chaos theory, and its applications in modern science. I look forward to future collaborations and the continued pursuit of scientific exploration.

Thank you, Chaima Zeraoulia, for your remarkable contributions and dedication to this research.

\end{document}